\documentclass[review]{elsarticle}

\usepackage{hyperref}

\journal{Journal of \LaTeX\ Templates}

\usepackage[mathscr]{euscript}
\usepackage{amsmath, amsthm, amssymb}

\usepackage{mathtools}

\usepackage{xcolor}

\usepackage{stmaryrd}  %

\newcommand{\paren}[1]{\left( #1 \right)}

\newcommand{\set}[1]{\left\{ #1 \right\}}

\newcommand{\beq}{\begin{eqnarray*}}
\newcommand{\eeq}{\end{eqnarray*}}
\newcommand{\beqn}{\begin{eqnarray}}
\newcommand{\eeqn}{\end{eqnarray}}
\newcommand{\ben}{\begin{enumerate}}
\newcommand{\een}{\end{enumerate}}
\newcommand{\bit}{\begin{itemize}}
\newcommand{\eit}{\end{itemize}}
\providecommand{\hide}[1]{}

\newcommand{\inv}{^{-1}} %

\newcommand{\X}{\Omega}

\newcommand{\ddim}{\operatorname{ddim}}
\newcommand{\diam}{\operatorname{diam}}

\renewcommand{\phi}{\varphi}
\newcommand{\mst}{\operatorname{MST}}

\newtheorem{theorem}{Theorem}[section]

\begin{document}

\begin{frontmatter}

\title{Non-uniform packings}

\author[mymainaddress]{Lee-Ad Gottlieb}
\ead{leead@ariel.ac.il}

\author[mysecondaryaddress]{Aryeh Kontorovich\corref{mycorrespondingauthor}}
\cortext[mycorrespondingauthor]{Corresponding author}
\ead{karyeh@cs.bgu.ac.il}

\address[mymainaddress]{Ariel University, Ariel, Israel}
\address[mysecondaryaddress]{Ben-Gurion University of the Negev, Beer Sheva, Israel}

\begin{abstract}
We generalize the classical notion of packing a set by balls
with identical radii to the case where the radii may be different.
The largest number of such balls that fit inside the set
without overlapping is called its {\em non-uniform packing number}.
We show that the non-uniform packing number can be upper-bounded
in terms of the {\em average} radius of the balls,
resulting in bounds of the familiar classical form.
\end{abstract}


\end{frontmatter}

\section{Introduction}
Packing numbers (along with their dual notion of covering numbers)
provide a quantitative notion of compactness for a totally bounded
metric space and make a pervasive appearance in 
empirical processes
\cite{MR3837109},
learning theory \cite{hksw19},
and information theory
\cite{lint_introduction_1999}, among other fundamental results.
We note in passing that
violating the triangle inequality
destroys the covering-packing duality,
and packing numbers emerge as 
the more fundamental notion,
at least in
a learning-theoretic setting
\cite{gkn-jmlr17+aistats}.

We refer the reader to \cite{MR3728284} for basic metric-space
notions such as total boundedness and compactness.
Briefly,
a {\em metric space}
$(\X,\rho)$ is a set endowed with a positive symmetric
function, which additionally satisfies the triangle inequality. For $r>0$, a set $A\subseteq\X$ is said to be $r$-separated if $\rho(a,a')>r$ for all distinct $a,a'\in A$.
The $r$-packing number of $\X$, which we denote by
$M(r)$, is the maximum cardinality of any $r$-separated subset of
$\X$ (and is finite whenever $\X$ is totally bounded).

We will also need the notion of the {\em doubling dimension}
of a metric space; the latter is known to be of critical
algorithmic 
\cite{KL04, HM06, CG06, DBLP:conf/stoc/Talwar04, DBLP:journals/siamcomp/BartalGK16}
and learning-theoretic importance
\cite{Bshouty2009323,
DBLP:journals/tit/GottliebKK14+colt,
KpotufeDasgupta2012,GottliebKK13tcs+alt,
gkn-ieee18+nips}.
Denote by $B(x,r)=\set{x'\in\X:\rho(x,x')\le r}$
the (closed) $r$-ball about $x$. If
there is a $D<\infty$ such that
every $r$-ball in $\X$
is contained in the union of some $D$ $r/2$-balls,
the metric space $(\X,\rho)$ is said to be {\em doubling}.
Its {\em doubling dimension} is defined as $\ddim(\X)=\ddim(\X,\rho)=:\log_2D^*$,
where $D^*$ is the smallest $D$ verifying the doubling property.
It is well-known \citep{KL04,GottliebKK13tcs+alt} that
\beqn
\label{eq:ddim-pack}
M(r)\le\paren{\frac{2\diam(\X)}{r}}^{\ddim(\X)},
\qquad r>0
,
\eeqn
where $\diam(\X)=\sup_{x,x'\in\X}\rho(x,x')$.
Further, (\ref{eq:ddim-pack}) is tight,
as witnessed by the example of $n$ equidistant points,
with $r$ as $(1-\varepsilon)$ times their common distance, for $\varepsilon$ arbitrarily small; in this case, $\ddim(\X) =\log_2 n$.

We now refine the notion of $r$-separated
sets to take the individual
inter-point distances into account.
For 
$A\subseteq\X$
and
$R:A\to (0,\infty)$, we say that
$A$ is $R$-separated if
for all $a\in A$,
\beqn
\inf_{a'\in A\setminus\set{a}}\rho(a,a')
> R(a).
\eeqn
In words, for each $a\in A$, its closest neighbor in $A$
is at least $R(a)$-away. The uniform special case $R(a)\equiv r$
recovers the classical notion of $r$-separation.

We are now ready to state our main result:
\begin{theorem}
\label{thm:main}
If $(\X,\rho)$ is a doubling space
and $A\subseteq\X$ is finite and $R$-separated, then
\beq
|A| 
\le \paren{\frac{5\diam(A)}{\bar r}}^{\min \{ \ddim(A), \ddim(\X) \}},
\eeq
where $\bar r:=|A|\inv\sum_{a\in A}R(a)$
is the average separation radius.
\end{theorem}

Observe that for the uniform special case
$R(a)\equiv r$,
Theorem~\ref{thm:main} recovers
(\ref{eq:ddim-pack}) up to constants.
We note that while $\ddim(A)$ may be arbitrarily 
smaller than $\ddim(\X)$, it may also be larger,
as $A$ may lack points used as ball centers in 
coverings of $\X$.
However, \cite{GK-13}
demonstrated that
for all $A\subseteq\X$,
we have
$\ddim(A)\le2\ddim(\X)$.

\paragraph{Related work}
The only tangentially relevant works we found
study the algorithmic 
\cite{MR3548502,MR3723337}
and game-theoretic \cite{MR3162285}
aspects
of optimization problems
involving packing
different-sized items
under various bin
constraints.
The results proved here were
early precursors to attempts
at defining a useful notion of
average Lipschitz smoothness,
but that line of research
ended up using entirely unrelated
techniques \cite{DBLP:AGK-20}.

\section{Proofs}

Before proving Theorem~\ref{thm:main} in its full
generality, we find it instructive to prove
the special case
where $(\X,\rho)$
is the unit ball of a $d$-dimensional normed space.
Any such space can be endowed with the Lebesgue measure
$\mu$ such that the $\mu$-volume of any $r$-ball is
$Cr^d$, where $C$ depends on 
the norm
and $d$ only.
Now if $A\subset\X$ is $R$-separated, then the balls
$B(a,R(a)/2)$ are all disjoint and contained in
$B(0,2)$.
Thus, the total volume of these balls is at most
$C2^d$ and at least
\beq
C\sum_{a\in A}(R(a)/2)^d.
\eeq
Combining these, we get the inequality
\beq
\sum_{a\in A}R(a)^d \le 4^d.
\eeq
Jensen's inequality implies that
\beq
\bar r^d=
\paren{ |A|\inv\sum_{a\in A}R(a)}^d
\le
|A|\inv\sum_{a\in A}R(a)^d,
\eeq
whence
\beq
\bar r^d
\le
|A|\inv 4^d.
\eeq
Solving for $|A|$ yields the bound
\beqn
\label{eq:volumetric}
|A|\le (4/\bar r)^d,
\eeqn
which recovers, up to constants, the classic
volumetric packing bounds
(see, e.g., \cite[Lemma 5.7]{9781108498029})
in the uniform special case $R(a)\equiv r$.
The aforementioned lemma
shows that
$d$-dimensional
normed spaces
have $\ddim\le d\log_2 6$.

Although the bound (\ref{eq:volumetric}) is very much
in the spirit of Theorem~\ref{thm:main}, the volumetric
technique does not extend to general metric spaces.
We will instead make use of weighted spanning tress.

\begin{proof}[Proof of Theorem~\ref{thm:main}]
There is no loss of generality in
normalizing all of the distances so that
$\diam(A)=1$.
Put $N:=|A|$ and
$\bar r:=N\inv\sum_{a\in A}R(a)$.
We will show that
\beqn
\label{eq:N-ub}
N < (5/\bar r)^{\min \{ \ddim(A), \ddim(\X) \}},
\eeqn
which proves the Theorem statement.

To prove (\ref{eq:N-ub}), let
the Minimum Spanning Tree
of $A$,
denoted
$\mst(A)$,
be rooted at a point $t \in A$
for which $R(t)$ is minimal, and it must be that 
$R(t) \le 1$. 
Let $E$ be the edge-set of $\mst(A)$, and denote
the length of each edge $e \in E$ by $l(e)$.
Further define 
$l(E) = \sum_{e \in E} l(e)$.
Now assign each edge of $E$ to the endpoint farthest from the root $t$;
this assigns a single edge to each point in the tree,
except to the root $t$. 
Let the edge assigned to a point $a \in A$
be $e(a)$, and for convenience we will say that
$e(t)$ is an edge of infinite length.
Trivially, the edge assigned to
each endpoint cannot be shorter than the distance from the endpoint to its
nearest neighbor in $A$, so
$R(a) \le l(e(a))$ for all $a \in A$.
It follows that
$N\bar r 
= \sum_{a \in A}R(a)
= \sum_{a \ne t \in A}R(a) + R(t)
\le l(E) + 1$.

Now Talwar
\cite[Lemma 6]{DBLP:conf/stoc/Talwar04}
(see also
\cite[Proposition 12]{DBLP:journals/jacm/Arora98}) has shown that the length of the MST on any set $A \in \X$ 
of $N$ points is at most
$$4\diam(A) N ^{1-1/\min \{ \ddim(A), \ddim(\X) \} }.$$
As we have taken the diameter to be bounded by $1$, we have
$\bar r \le (l(E) + 1)/N
    \le 4N ^{-1/\min \{ \ddim(A), \ddim(\X) \}} + 1/N
    <   5 N ^{-1/\min \{ \ddim(A), \ddim(\X) \}}$.
The bound claimed in (\ref{eq:N-ub}) follows.

\end{proof}

\bibliography{refs}

\begin{thebibliography}{10}
\expandafter\ifx\csname url\endcsname\relax
  \def\url#1{\texttt{#1}}\fi
\expandafter\ifx\csname urlprefix\endcsname\relax\def\urlprefix{URL }\fi
\expandafter\ifx\csname href\endcsname\relax
  \def\href#1#2{#2} \def\path#1{#1}\fi

\bibitem{MR3837109}
R.~Vershynin, \href{https://doi.org/10.1017/9781108231596}{High-dimensional
  probability}, Vol.~47 of Cambridge Series in Statistical and Probabilistic
  Mathematics, Cambridge University Press, Cambridge, 2018, an introduction
  with applications in data science, With a foreword by Sara van de Geer.
\newblock \href {http://dx.doi.org/10.1017/9781108231596}
  {\path{doi:10.1017/9781108231596}}.
\newline\urlprefix\url{https://doi.org/10.1017/9781108231596}

\bibitem{hksw19}
S.~Hanneke, A.~Kontorovich, S.~Sabato, R.~Weiss,
  \href{http://arxiv.org/abs/1906.09855}{Universal bayes consistency in metric
  spaces} (2019).
\newblock \href {http://arxiv.org/abs/1906.09855} {\path{arXiv:1906.09855}}.
\newline\urlprefix\url{http://arxiv.org/abs/1906.09855}

\bibitem{lint_introduction_1999}
J.~H. Lint,
  \href{http://public.ebookcentral.proquest.com/choice/publicfullrecord.aspx?p=3092660}{Introduction
  to {Coding} {Theory}}, Springer Berlin Heidelberg : Imprint : Springer,
  Berlin, Heidelberg, 1999, oCLC: 840292572.
\newline\urlprefix\url{http://public.ebookcentral.proquest.com/choice/publicfullrecord.aspx?p=3092660}

\bibitem{gkn-jmlr17+aistats}
L.-A. Gottlieb, A.~Kontorovich, P.~Nisnevitch, Nearly optimal classification
  for semimetrics (extended abstract: {AISTATS} 2016), Journal of Machine
  Learning Research.

\bibitem{MR3728284}
J.~R. Munkres, Topology, Prentice Hall, Inc., Upper Saddle River, NJ, 2000,
  second edition of [ MR0464128].

\bibitem{KL04}
R.~Krauthgamer, J.~R. Lee, Navigating nets: {S}imple algorithms for proximity
  search, in: 15th Annual ACM-SIAM Symposium on Discrete Algorithms, 2004, pp.
  791--801.

\bibitem{HM06}
S.~{Har-Peled}, M.~Mendel, \href{http://link.aip.org/link/?SMJ/35/1148/1}{Fast
  construction of nets in low-dimensional metrics and their applications}, SIAM
  Journal on Computing 35~(5) (2006) 1148--1184.
\newblock \href {http://dx.doi.org/10.1137/S0097539704446281}
  {\path{doi:10.1137/S0097539704446281}}.
\newline\urlprefix\url{http://link.aip.org/link/?SMJ/35/1148/1}

\bibitem{CG06}
R.~Cole, L.-A. Gottlieb, Searching dynamic point sets in spaces with bounded
  doubling dimension, in: STOC, 2006, pp. 574--583.

\bibitem{DBLP:conf/stoc/Talwar04}
K.~Talwar, \href{http://doi.acm.org/10.1145/1007352.1007399}{Bypassing the
  embedding: algorithms for low dimensional metrics}, in: Proceedings of the
  36th Annual {ACM} Symposium on Theory of Computing, Chicago, IL, USA, June
  13-16, 2004, 2004, pp. 281--290.
\newblock \href {http://dx.doi.org/10.1145/1007352.1007399}
  {\path{doi:10.1145/1007352.1007399}}.
\newline\urlprefix\url{http://doi.acm.org/10.1145/1007352.1007399}

\bibitem{DBLP:journals/siamcomp/BartalGK16}
Y.~Bartal, L.~Gottlieb, R.~Krauthgamer,
  \href{https://doi.org/10.1137/130913328}{The traveling salesman problem:
  Low-dimensionality implies a polynomial time approximation scheme}, {SIAM} J.
  Comput. 45~(4) (2016) 1563--1581.
\newblock \href {http://dx.doi.org/10.1137/130913328}
  {\path{doi:10.1137/130913328}}.
\newline\urlprefix\url{https://doi.org/10.1137/130913328}

\bibitem{Bshouty2009323}
N.~H. Bshouty, Y.~Li, P.~M. Long,
  \href{http://www.sciencedirect.com/science/article/B6WJ0-4VH4DPR-1/2/ec8c49f50cae69c0e92f71f7a4be6691}{Using
  the doubling dimension to analyze the generalization of learning algorithms},
  Journal of Computer and System Sciences 75~(6) (2009) 323 -- 335.
\newblock \href {http://dx.doi.org/DOI: 10.1016/j.jcss.2009.01.003}
  {\path{doi:DOI: 10.1016/j.jcss.2009.01.003}}.
\newline\urlprefix\url{http://www.sciencedirect.com/science/article/B6WJ0-4VH4DPR-1/2/ec8c49f50cae69c0e92f71f7a4be6691}

\bibitem{DBLP:journals/tit/GottliebKK14+colt}
L.~Gottlieb, A.~Kontorovich, R.~Krauthgamer,
  \href{http://dx.doi.org/10.1109/TIT.2014.2339840}{Efficient classification
  for metric data (extended abstract: {COLT} 2010)}, {IEEE} Transactions on
  Information Theory 60~(9) (2014) 5750--5759.
\newblock \href {http://dx.doi.org/10.1109/TIT.2014.2339840}
  {\path{doi:10.1109/TIT.2014.2339840}}.
\newline\urlprefix\url{http://dx.doi.org/10.1109/TIT.2014.2339840}

\bibitem{KpotufeDasgupta2012}
S.~Kpotufe, S.~Dasgupta, \href{http://dx.doi.org/10.1016/j.jcss.2012.01.002}{A
  tree-based regressor that adapts to intrinsic dimension}, J. Comput. Syst.
  Sci. 78~(5) (2012) 1496--1515.
\newblock \href {http://dx.doi.org/10.1016/j.jcss.2012.01.002}
  {\path{doi:10.1016/j.jcss.2012.01.002}}.
\newline\urlprefix\url{http://dx.doi.org/10.1016/j.jcss.2012.01.002}

\bibitem{GottliebKK13tcs+alt}
L.-A. Gottlieb, A.~Kontorovich, R.~Krauthgamer, Adaptive metric dimensionality
  reduction (extended abstract: {ALT} 2013), Theoretical Computer Science
  (2016) 105--118.

\bibitem{gkn-ieee18+nips}
L.~Gottlieb, A.~Kontorovich, P.~Nisnevitch,
  \href{https://doi.org/10.1109/TIT.2018.2822267}{Near-optimal sample
  compression for nearest neighbors (extended abstract: {NIPS} 2014)}, {IEEE}
  Trans. Information Theory 64~(6) (2018) 4120--4128.
\newblock \href {http://dx.doi.org/10.1109/TIT.2018.2822267}
  {\path{doi:10.1109/TIT.2018.2822267}}.
\newline\urlprefix\url{https://doi.org/10.1109/TIT.2018.2822267}

\bibitem{GK-13}
L.-A. Gottlieb, R.~Krauthgamer, Proximity algorithms for nearly doubling
  spaces, SIAM J. Discrete Math. 27~(4) (2013) 1759--1769.

\bibitem{MR3548502}
Y.~G. Stoyan, G.~Sha\u{\i}tkhauer, G.~N. Yas'kov,
  \href{https://doi.org/10.1007/s10559-016-9842-1}{Packing unequal spheres into
  different containers}, Kibernet. Sistem. Anal. 52~(3) (2016) 97--105.
\newblock \href {http://dx.doi.org/10.1007/s10559-016-9842-1}
  {\path{doi:10.1007/s10559-016-9842-1}}.
\newline\urlprefix\url{https://doi.org/10.1007/s10559-016-9842-1}

\bibitem{MR3723337}
A.~Ene, S.~Har-Peled, B.~Raichel,
  \href{https://doi.org/10.1137/120898413}{Geometric packing under nonuniform
  constraints}, SIAM J. Comput. 46~(6) (2017) 1745--1784.
\newblock \href {http://dx.doi.org/10.1137/120898413}
  {\path{doi:10.1137/120898413}}.
\newline\urlprefix\url{https://doi.org/10.1137/120898413}

\bibitem{MR3162285}
W.~Kern, X.~Qiu, \href{https://doi.org/10.1016/j.dam.2012.08.002}{Note on
  non-uniform bin packing games}, Discrete Appl. Math. 165 (2014) 175--184.
\newblock \href {http://dx.doi.org/10.1016/j.dam.2012.08.002}
  {\path{doi:10.1016/j.dam.2012.08.002}}.
\newline\urlprefix\url{https://doi.org/10.1016/j.dam.2012.08.002}

\bibitem{DBLP:AGK-20}
Y.~Ashlagi, L.~Gottlieb, A.~Kontorovich,
  \href{https://arxiv.org/abs/2007.06283}{Functions with average smoothness:
  structure, algorithms, and learning} (2020).
\newblock \href {http://arxiv.org/abs/2007.06283} {\path{arXiv:2007.06283}}.
\newline\urlprefix\url{https://arxiv.org/abs/2007.06283}

\bibitem{9781108498029}
M.~J. Wainwright,
  \href{https://www.amazon.com/High-Dimensional-Statistics-Non-Asymptotic-Statistical-Probabilistic/dp/1108498027?SubscriptionId=AKIAIOBINVZYXZQZ2U3A&tag=chimbori05-20&linkCode=xm2&camp=2025&creative=165953&creativeASIN=1108498027}{High-Dimensional
  Statistics: A Non-Asymptotic Viewpoint (Cambridge Series in Statistical and
  Probabilistic Mathematics)}, Cambridge University Press, 2019.
\newline\urlprefix\url{https://www.amazon.com/High-Dimensional-Statistics-Non-Asymptotic-Statistical-Probabilistic/dp/1108498027?SubscriptionId=AKIAIOBINVZYXZQZ2U3A&tag=chimbori05-20&linkCode=xm2&camp=2025&creative=165953&creativeASIN=1108498027}

\bibitem{DBLP:journals/jacm/Arora98}
S.~Arora, \href{http://doi.acm.org/10.1145/290179.290180}{Polynomial time
  approximation schemes for euclidean traveling salesman and other geometric
  problems}, J. {ACM} 45~(5) (1998) 753--782.
\newblock \href {http://dx.doi.org/10.1145/290179.290180}
  {\path{doi:10.1145/290179.290180}}.
\newline\urlprefix\url{http://doi.acm.org/10.1145/290179.290180}

\end{thebibliography}

\end{document}